\input amstex
\documentstyle{amsppt}
\def \al{\alpha}
\def \bt{\beta}
\def \gm{\gamma}

\def \sd #1#2{{#1}^{\scriptstyle{(#2)}}}

\def \QQ{\Bbb {Q}}
\def \NN{\Bbb {N}}
\def \ZZ{\Bbb {Z}}

\def \mm{\Cal M_m^+}
\def \mn{\Cal M_m}

\pagewidth{125mm} 
\pageheight{185mm} 
\parindent=8mm 
\frenchspacing 
\NoBlackBoxes

\topmatter
\title The ring of multisymmetric functions
\endtitle
\author Francesco Vaccarino\endauthor

\address Dipartimento di  Matematica - Politecnico di Torino
- Corso Duca degli Abruzzi 24 - 10129 - Torino - Italy\endaddress 

\email vaccarino\@syzygie.it \endemail



\abstract 
Let $R$ be a commutative ring and let $n,m$ be two positive
integers. Let $A_R(n,m)$ be the polynomial ring in the commuting
independent variables $ x_i(j)$ with $i=1, \dots ,m\,;j=1, \dots  ,n$
and coefficients in $R$.  The symmetric group on $n$ letters $S_n$
acts on $A_R(n,m)$ by means of $\sigma(x_i(j))=x_i(\sigma(j))$ for all
$\sigma\in S_n$ and $i=1, \dots ,m\,;j=1, \dots  ,n$. Let us denote by
$A_R(n,m)^{S_n}$ the ring of invariants for this action: its elements
are usually called multisymmetric functions and they are the usual
symmetric functions when $m=1$. In this paper we give a
presentation of $A_R(n,m)^{S_n}$ in terms of generators and relations
that holds for any $R$ and any $n,m$, thereby answering a classical
question. 
\endabstract
\endtopmatter

\document

\head Introduction \endhead
Let $R$ be a commutative ring and let $n,m$ be two positive
integers. Let $A_R(n,m)$ be the polynomial ring in the commuting
independent variables $ x_i(j)$ with $i=1, \dots ,m\,;j=1, \dots  ,n$
and coefficients in $R$.  The symmetric group on $n$ letters $S_n$
acts on $A_R(n,m)$ by means of $\sigma(x_i(j))=x_i(\sigma(j))$ for all
$\sigma\in S_n$ and $i=1, \dots ,m\,;j=1, \dots  ,n$. Let us denote by
$A_R(n,m)^{S_n}$ the rings of invariants for this action: its elements
are usually called multisymmetric functions and they are the usual
symmetric functions when $m=1$. In this case, 
$A_R(n,1)\cong R[x_1,x_2,\dots ,x_n]$, and 
$R[x_1,x_2,\dots ,x_n]^{S_n}$ is freely generated by the elementary
symmetric functions $e_1,\dots,e_n$ given by the equality 
$$
\sum_{k=0}^n t^k e_k:=\prod_{i=1}^n(1+tx_i).\eqno(0.1)
$$
Here $e_0=1$ and $t$ is a commuting independent variable (see {\cite{M}}).
Furthermore one has
$$e_k(x_1,\dots,x_n)=\sum_{i_1<i_2<\dots <i_k\leq n}x_{i_1}x_{i_2}\cdots
x_{i_k}\eqno(0.2)$$

\noindent
Unless otherwise stated, we now assume that $m>1$. We first obtain
generators of the ring $A_R(n,m)^{S_n}$. 

\noindent
Let $A_R(m):=R[y_1,\dots,y_m]$, where $y_1,\dots,y_m$ are commuting
independent variables, let $f=f(y_1,\dots,y_m)\in A_R(m)$ and define
$$
f(j):=f(x_1(j),\dots ,x_m(j))\,{\text{ for }}\, 1 \leq j\leq n.
\eqno(0.3)
$$  
Notice that $f(j)\in A_R(n,m)$ for all $1 \leq j\leq n $ and that
$\sigma(f(j))=f(\sigma(j))$, for all $\sigma\in S_n$ and
$j=1,\dots,n$.  

\noindent
Define $e_k(f):=e_k(f(1),f(2),\dots,f(n))$ i.e.
$$\sum_{k=0}^n t^k e_k(f):=\prod_{i=1}^n(1+tf(i)),\eqno(0.4)
$$
where $t$ is a commuting independent variable. Then 
$e_k(f)\in A_R(n,m)^{S_n}$. 

\noindent
One may think about the $y_i$ as diagonal matrices in the following
sense: let $M_n(A_R(n,m))$ be the full ring of $n\times n$ matrices
with coefficients in $A_R(n,m)$. Then there is an embedding  
$$\rho_n:A_R(m)\hookrightarrow M_n(A_R(n,m))\eqno(0.5)$$ 
given by 
$$\rho_n(y_i):=
\pmatrix x_i(1)&0&\dots&0\\
 0&x_i(2)&\dots&0\\
	0&0&\dots&x_i(n)
\endpmatrix {\text{ for }} i=1,\dots,m.\eqno(0.6)
$$
Now $(0.4)$ gives
$$
\sum_{k=0}^n t^k e_k(f) =
\prod_{j=1}^n(1+t\rho_n(f)_{jj})=det(1+t\rho_n(f)),\eqno(0.7)
$$ 
where $det(-)$ is the usual determinant of $n\times n$ matrices. 

\noindent
Let $\Cal M_m$ be the set of monomials in $A_R(m)$. For 
$\mu\in\Cal M_m$ let $\partial_i(\mu)$ denote the degree of $\mu$ in
$y_i$, for all $i=1,\dots ,m$. We set 
$$
\partial(\mu):=(\partial_1(\mu),\dots ,\partial_m(\mu))\eqno(0.8)
$$  
for its multidegree. The total degree of $\mu$ is $\sum_i\partial_i(\mu)$.
Let $\mm$ be the set of monomials of positive degree.  
A monomial $\mu\in \mm$ is called {\it{primitive}} it is not a power
of another one. We denote by $\frak M_m^+$ the set of primitive monomials.
We define an $S_n$ invariant multidegree on $A_R(n,m)$ by setting
$\partial(x_i(j))=\partial(y_i)\in \NN^m$ for all $1\leq j\leq n$ and
$1\leq i\leq m$. If $f\in A_R(m)$ is homogeneous of total degree $l$,
then $e_k(f)$ has total degree $kl$ (for all $k$ and $n$). 

\noindent
We are now in a position to state the first part of our result (recall
that $m>1$). 

\proclaim{Theorem 1 (generators)}
The ring of multisymmetric functions $A_R(n,m)^{S_n}$ is generated by
the $e_k(\mu)$, where $\mu\in \frak M_m^+$, $k=1,\dots n$ and the total
degree of $e_k(\mu)$ is less or equal than $n(m-1)$. 
If $n=p^s$ is a power of a prime and $R=\ZZ$ or $p\cdot 1_R=0$, then at
least one generator has degree equal to $n(m-1)$.  

\noindent
If $R\supset \QQ$ then $A_R(n,m)^{S_n}$ is generated by the $e_1(\mu)$,
where $\mu\in \mm$ and the degree of $\mu$ is less or equal than $n$.  
\endproclaim

\noindent
To obtain the relations between these generators, we need more
notation on (multi)symmetric functions. 

\noindent
The action of $S_n$ on $A_R(n,1)\cong R[x_1,x_2,\dots ,x_n]$ preserves
the usual degree. We denote by $\Lambda_{R,n}^k$ the $R$-submodule of 
invariants of degree $k$.

\noindent
Let $q_n: R[x_1,x_2,\dots ,x_n]@>>> R[x_1,x_2,\dots ,x_{n-1}]$ be
given by $x_n\mapsto 0$ and $x_i\mapsto x_i$, for
$i=1,\dots,n-1$. This map sends  $\Lambda_{n,R}^k$ to
$\Lambda_{n-1,R}^k$ and it is easy to see that
$\Lambda_{n,R}^k\cong\Lambda_{k,R}^k$ for all $n \geq k$. Denote by
$\Lambda^k_R$ the limit of the inverse system obtained in this way.  

\noindent
The ring $\Lambda_R:=\bigoplus_{k\geq 0}\Lambda^k_R$ is called the
ring of {\it{symmetric functions}} (over $R$). 

\noindent
It can be shown {\cite{M}} that $\Lambda_R$ is a polynomial ring,
freely generated by the (limits of the) $e_k$, that are given by 
$$
\sum_{k=0}^{\infty} t^k e_k:=\prod_{i=1}^{\infty}(1+tx_i). \eqno(0.9)
$$
Furthermore the kernel of the natural projection 
$\pi_n:\Lambda_R @>>> \Lambda_{n,R}$ is generated by the $e_{n+k}$,
where $k\geq 1$. 
 
\noindent
In a similar way we build a limit of multisymmetric functions. 
For any $a\in \NN^m$ we set $A_R(n,m,a)$ for the linear span of the
monomials of multidegree $a$. 
One has 
$$
A_R(n,m)=\bigoplus_{a\in \NN^m}A_R(n,m,a).\eqno(0.10)
$$ 

\noindent
Let $\pi_{n}:A_R(n,m)@>>> A_R(n-1,m)$ be given by
$$\pi_{n}(x_i(j))=\cases 0 \;&{\text{if}}\,\, j=n \\
x_i(j) \;&{\text{if}}\,\, j\leq n-1
\endcases \,\,\,\,\,\;\;\;\;\;\;\;{\text{ for all }}\, i.\eqno(0.11)
$$
Then (see (3.5)) we prove that, for all $a\in\NN^m$
$$ \pi_{n}(A_R(n,m,a)^{S_n})=A_R(n-1,m,a)^{S_{n-1}}.\eqno(0.12)$$
\noindent
For any $a \in \NN^m$ set 
$$A_R(\infty,m,a):=\lim_{\leftarrow } A_R(n,m,a)^{S_n}, \eqno(0.13)$$ 
where the projective limit is taken with respect to $n$ over the
projective system  
$(A_R(n,m,a)^{S_n},\pi_{n})$.

\noindent
Set 
$$A_R(\infty,m):=\bigoplus_{a\in \NN^m} A_R(\infty,m,a).\eqno(0.14)$$

\noindent
We set, by abuse of notation, $$e_k(f):=\lim_{\leftarrow }e_k(f)\in A_R(\infty,m) \eqno(0.15)$$ 
with $k\in \NN$ and $f\in A(m)^+$, the augmentation ideal, i.e.
$$
\sum_{k=0}^{\infty} t^k e_k(f):=\prod_{j=1}^{\infty}(1+tf(j)). \eqno(0.16)
$$

\noindent
Then $e_k$ is a homogeneous polynomial of degree $k$. Now, if
$f=\sum_{\mu\in\mm}\lambda_{\mu}\mu$, we set
$$e_k(f):=\sum_{\al}\lambda^{\al}e_{\al}\eqno(0.16)$$
where $\al:=(\al_{\mu})_{\mu\in\mm}$ is such that $\al_{\mu}\in\NN$, 
$\sum_{\mu\in\mm}\al_{\mu}\leq k$ and
$\lambda^{\al}:=\prod_{\mu\in\mm}\lambda^{\al_{\mu}}$. 

\noindent
We can now state the second part of our main result.

\proclaim{Theorem 2 (relations)}
\roster
\item
The ring $A_R(\infty,m)$ is a polynomial ring, freely
generated by the (limits of) the $e_k(\mu)$, where 
$\mu\in\frak M_m^+$ and $k\in \NN$.  

\noindent
The kernel of the natural projection 
$$A_R(\infty,m)@>>> A_R(n,m)^{S_n}$$
is generated as $R$-module by the coefficients $e_{\al}$ of the
elements 
$$
e_{n+k}(f), \,{\text{ where }} \, k\geq 1 \,{\text{ and }}\, f\in A_R(m)^+.
$$ 
\item If $R\supset \QQ$ then $A_R(\infty,m)$ is freely generated by
the $e_1(\mu)$, where $\mu\in\mm$. 

\noindent
The kernel of the natural projection is generated as an ideal by the
$e_{n+1}(f)$, where $f\in A_R(m)^+$. 
\endroster
\endproclaim

\noindent
In Dalbec's paper {\cite{D}} generators and relations are found in the
case where $R \supset \QQ$. The relations found there are actually the
same we find: indeed what Dalbec calls {\it{monomial multisymmetric
functions}} are exactly those $e_{\al}$ we introduced in (0.17), so
that his Proposition 1.9 is a special case of our Proposition 3.1(1)
when $R \supset \QQ$. Another paper on this theme, giving a minimal
presentation when the base ring is a characteristic 2 field, is
{\cite{A}}. Again, its main results on multisymmetric functions are a
corollary of ours when $R$ is a characteristic 2 field. 

\noindent
The results of this paper were presented in 1997 at a congress on
algebraic groups representations in Ascona (CH) organized by
H.P.Kraft. They are published only now for personal reasons.

\head 1. Notations and basic facts\endhead

\noindent
The monomials of $A_R(n,m)$ form a $R$-basis, permuted by the action
of $S_n$. Thus, the sums of monomials over the orbits form a $R$-basis
of the ring of multisymmetric functions. We now introduce some
notation and preliminary results concerning these functions and orbit sums.

Let $k\in \NN$, we denote by $\bold{f}$ the sequence 
$(f_1\dots,f_k)$ in $A_R(m)$ and by $\al$ the element
$(\al_1,\dots,\al_k)\in\NN^k$, where $\sum\al_j \leq n$ .  
Let $t_1,\dots,t_k$ be commuting independent variables, we set as
usual $t^{\al}:=\prod_{i}t_i^{\al_i}$. We define elements
$e_{\al}(\bold{f})\in A_R(n,m)^{S_n}$ by 
$$
\sum_{\al}t^{\al}e_{\al}(\bold{f}) :=
det(1+\sum_ht_h\rho_n(f_h))=\prod_{i=1}^n(1+\sum_ht_h f_h(i)).
\eqno(1.1)
$$

\example{Example 1.1}
Let $n=3$ and $f,g\in A_R(m)$ then 
$$e_{(2,1)}(f,g)= f(1)f(2)g(3)+f(1)g(2)f(3)+g(1)f(2)f(3).$$
If $n=4$ then
$$\aligned e_{(2,1)}(f,g) = &f(1)f(2)g(3)+f(1)g(2)f(3)+g(1)f(2)f(3)+\\
              &f(1)f(2)g(4)+f(1)g(2)f(4)+g(1)f(2)f(4)+\\
              &f(1)f(3)g(4)+f(1)g(3)f(4)+g(1)f(3)f(4)+\\
                 &f(2)f(3)g(4)+f(2)g(3)f(4)+g(2)f(3)f(4)\endaligned$$
                 \endexample

\noindent
Let $k=m$ and $f_j=y_j$ for $j=1,\dots,m$, then the
$e_{\al}(\bold{y})=e_{(\al_1,\dots,\al_m)}(y_1,\dots,y_m)$  
where $\sum\al_j\leq n$ are the well--known {\it{elementary
multisymmetric functions}}. These generate $A_R(n,m)^{S_n}$ when
$R\supset\QQ$ (see {\cite{G}} or {\cite{W}}), and satisfy 
$$
\sum_{\al}t^{\al}e_{\al}(\bold{y})
=det(1+\sum_jt_j\rho_n(y_j))
=\prod_{i=1}^n(1+\sum_{j=1}^m t_jx_j(i)).\eqno(1.2)
$$

\noindent
\proclaim{Lemma 1.2}
The multisymmetric function $e_{(\al_1,\dots,\al_k)}(f_1,\dots,f_k)$ is the orbit sum (under the considered action of $S_n$) of $$f_1(1)f_1(2)\cdots f_1(\al_1)f_2(\al_1+1)\cdots f_2(\al_1+\al_2)\cdots f_k(\sum_h \al_h).$$
\endproclaim
\demo{Proof}
Let $E$ be the set of mappings 
$\phi:\left\{1,\dots,n\right\}\rightarrow\left\{1,\dots,k+1\right\}$.
We define a mapping $\phi\mapsto\phi^*$ of $E$ into $\NN^{k+1}$ by
putting $\phi^*(i)$ equal to the cardinality of $\phi^{-1}(i)$. For
two elements $\phi_1,\phi_2$ of $E$, to satisfy $\phi_1^*=\phi_2^*$ it
is necessary and sufficient that there should exist $\sigma\in S_n$
such that $\phi_2=\phi_1\circ\sigma$. Set $f_{k+1}:=1_R$ and
$E(\al):=\left\{\phi\in E \mid \phi^*=(\al_1,\dots,\al_k,
n-\sum_i\al_i)\right\}$, then we have 
$$e_{\alpha}(\bold{f})
=\sum_{\phi\in E(\al)}f_{\phi(1)}(1)f_{\phi(2)}(2)\cdots f_{\phi(n)}(n)\eqno(1.3)
$$
and the lemma is proved.\qed
\enddemo

\noindent
It is clear that 
$e_{(\al_1,\dots,\al_k)}(f_1,\dots,f_k)
=e_{(\al_{\tau(1)},\dots,\al_{\tau(k)})}(f_{\tau(1)},\dots,f_{\tau(k)})$
for all $\tau\in S_k$. If two entries are equal, say $f_1=f_2$, then,
by $(1.1)$
$$e_{(\al_1,\dots,\al_k)}(f_1,\dots,f_k)
=\frac{(\al_1+\al_2)!}{\al_1!\al_2!}
e_{(\al_1+\al_2,\dots,\al_k)}(f_1,f_3\dots,f_k).\eqno(1.4)$$  

\noindent
Let $\sd {\NN} {\mn^+}$ be the set of functions $\mn^+ @>>> \NN$ with
finite support. We set 
$$\mid \al \mid:=\sum_{\mu\in \mm} \al(\mu)\eqno(1.5)$$

\noindent
Let $\al\in\sd {\NN} {\mn^+}$, then there exist $k\in \NN$ and
$\mu_1,\dots,\mu_k\in\mn^+$ such that $\al(\mu_i)=\al_i\neq 0$ for
$i=1,\dots,k$ and $\al(\mu)=0$ when $\mu\neq \mu_1,\dots,\mu_k$. We
set  
$$e_{\al}:=e_{(\al_1,\dots,\al_k)}(\mu_1,\dots,\mu_k),\eqno(1.6)$$
i.e. we substitute $(\mu_1,\dots,\mu_k)$ to variables in the
elementary multisymmetric function
$e_{(\al_1,\dots,\al_k)}(y_1,\dots,y_k)$. 

\noindent
Then 
$$\sum_{\mid\al\mid\leq n}t^{\al}e_{\al}
=\prod_{i=1}^n(1+\sum_{\mu\in\mn^+}t_{\mu}\mu(i)), \eqno(1.7)$$
where $t_{\mu}$ are commuting independent variables indexed by
monomials and
$$t^{\al}:=\prod_{\mu\in\mn^+}t_{\mu}^{\al(\mu)}\eqno(1.8)$$ 
for all $\al\in \sd {\NN} {\mn^+}$. 

\noindent
If $\al\in \sd {\NN} {\mn^+}$ is such that $\al(\mu)=k$ for some
$\mu\in\mn^+$ and $\al(\nu)=0$ for all $\nu\in\mn^+$ with $\nu\neq\mu$, 
we see that $e_{\al}=e_k(\mu)$, the $k$-th elementary symmetric
function evaluated at $(\mu(1),\mu(2),\dots,\mu(n))$.

\proclaim{Lemma 1.3}
Given a monomial $\mu\in A_R(n,m)$, there exist $\mu_1,\dots,\mu_n\in
A_R(m)$ such that $\mu=\mu_1(1)\cdots\mu_n(n)$.  
\endproclaim

\demo{Proof}
Let $\mu=\prod_{ij}x_i(j)^{a_{ij}}$ then $\mu_j=\prod_iy_i^{a_{ij}}$
for $j=1,\dots,n$. 
\qed
\enddemo

\proclaim{Proposition 1.4}
The set 
$$\Cal B_{n,m,R}:=\{e_{\al}\;:\;\mid \al\mid\leq n\}$$ 
is a $R$-basis of $A_R(n,m)^{S_n}$.

\noindent
The set 
$$\Cal B_{n,m,a,R}:=\{e_{\al}\;:\; \mid\al\mid\leq n {\text{ and
}}\partial(e_{\al})=a\}$$  
is a $R$-basis of $A_R(n,m,a)^{S_n}$, for all $a\in \NN^m$.
\endproclaim

\demo{Proof}
By Lemma 1.2 and (1.6), the $e_{\al}$ are a complete system of
representatives (for the action of $S_n$) of the orbit sums of the
products 
$$\{\mu_1(1)\mu_2(2)\cdots\mu_n(n)\,:\,\mu_i\in\mn \,,\,i=1,\dots,n\}.$$ 
So the first statement follows by Lemma 1.3.

\noindent
Notice that
$\partial(e_{\al})=\sum_{\mu\in\mn^+}\al_{\mu}\partial(\mu)$ to prove
the second statement. 
\qed
\enddemo

\head 2. Generators \endhead

\noindent
Let us calculate the product between two elements
$e_{\al},e_{\beta}\in \Cal B_{n,m,R}$ of the basis $\Cal B_{n,m,R}$.

\noindent
\proclaim{Theorem 2.1 - Product Formula}

\noindent
Let $k,h\in \NN$, $f_1\dots ,f_k,g_1,\dots , g_h \in A_R(m)$ and
$t_1,\dots,t_k, s_1,\dots ,s_h$ be commuting independent variables.
Set as in (1.1)
$$e_{\al}(\bold{f}):=e_{(\al_1,\dots,\al_k)}(f_1,\dots,f_k) 
{\text{and }}e_{\bt}(\bold{g}):=e_{(\bt_1,\dots,\bt_h)}(g_1,\dots,g_h).$$ 
Then
$$e_{\al}(\bold{f})e_{\bt}(\bold{g})=
\sum_{\gm} e_{\gm}(\bold{f},\bold{g},\bold{fg}),$$
where $\bold{fg} := 
(f_1g_1,f_1g_2,\dots,f_1g_h,f_2g_1,\dots,f_2g_h,\dots,f_kg_h)$ and 

\noindent
$
\gamma :=
(\gamma_{10},\dots,\gamma_{k0},\gamma_{01},\dots,\gamma_{0h},
\gamma_{11},\gamma_{12},\dots,\gamma_{kh})$
 are such that
$$
\cases \gamma_{ij}\in \NN \\
\mid \gamma \mid \leq n \\
\sum_{j=0}^h \gamma_{ij}=\al_i \;\;{\text{for}}\;\; i=1,\dots,k\\
\sum_{i=0}^k \gamma_{ij}=\bt_j \;\; {\text{for}}\;\; j=1,\dots,h.
\endcases
$$
\endproclaim

\demo{Proof}
The result follows from
$$\aligned
&(\sum_{\sum\al_j\leq
n}\prod_{j=1}^kt_j^{\al_j}e_{\al}(\bold{f}))(\sum_{\sum \bt_l\leq
n}\prod_{l=1}^h s_l^{\bt_l}e_{\bt}(\bold{g}))= \\ 
&(\sum_{\al}t^{\al}e_{\al}(\bold{f}))(\sum_{\bt}s^{\bt}e_{\bt}(\bold{g}))=\\ 
&\prod_{i=1}^n(1+\sum_{j=1}^k t_j f_j(i))\prod_{i=1}^n(1+\sum_{l=1}^h
s_l g_l(i))=\\  
&\prod_{i=1}^n(1+\sum_{j=1}^k t_j f_j(i)+\sum_{l=1}^h s_l
g_l(i)+\sum_{j,l} t_js_l f_j(i)g_l(i)). 
\endaligned$$
Introduce the new variables $u_{jl}$ with $j=1,\dots,k$ and $l=1,\dots,h$, then
$$\aligned
&\prod_{i=1}^n(1+\sum_{j=1}^k t_j f_j(i)+\sum_{l=1}^h s_l
g_l(i)+\sum_{j,l} t_js_l f_j(i)g_l(i))=\\ 
&\prod_{i=1}^n(1+\sum_{j=1}^k t_j f_j(i)+\sum_{l=1}^h s_l
g_l(i)+\sum_{j,l} u_{jl}(i)g_l(i))=\\ 
&\sum_{\gamma}v^{\gamma}e_{\gamma}(\bold{f},\bold{g},\bold{fg})
\endaligned
$$
where $v$ is the cumulative variable $t,s,u$.
Then substitute $u_{jl}=t_js_l$ to obtain
$$\aligned
&\sum_{\gamma}v^{\gamma}e_{\gamma}(\bold{f},\bold{g},\bold{fg})=\\
&\sum_{\gamma} (\prod_{a=1}^k t_a^{\gamma_{a0}}\prod_{b=1}^h s_b^{\gamma_{0b}}
\prod_{a=1}^k \prod_{b=1}^h (t_a
s_b)^{\gamma_{ab}}e_{\gamma}(\bold{f},\bold{g},\bold{fg})), 
\endaligned$$
where $\bold{fg}=(f_1g_1,f_1g_2,\dots,f_kg_1,\dots,f_kg_h)$ and $\gm$
satisfy the condition of the theorem.
\enddemo

\example{Example 2.2}
Let us calculate in $A_R(2,3)^{S_2}$
$$e_{(1,1)}(a,b)e_{2}(c)=\sum_{0\leq k,h \leq 1}
e_{(1-k,1-h,2-k-h,h,k)}(a,b,c,ac,bc)=e_{(1,1)}(ac,bc),$$  
since $1-k+1-h+2-k-h+h+k=4-k-h\leq 2$.
\endexample

\noindent
\proclaim{Corollary 2.3}
Let $k\in \NN$, $a_1,\dots ,a_k\in A_R(m)$, 
$\al=(\al_1,\dots,\al_k)\in \NN^k$ with $\sum\al_j\leq n$. Then
$e_{(\al_1,\dots,\al_k)}(a_1,\dots,a_k)$ belongs to the subring of
$A_R(n,m)^{S_n}$ generated by the $e_i(\mu)$, where $i=1,\dots,n$ and
$\mu$ is a monomial in the $a_1,\dots,a_k$. 
\endproclaim

\demo{Proof}
We prove the claim by induction on $\sum_j\al_j$ (notice that 
$1 \leq k\leq \sum_j\al_j$) assuming that $\al_i>0$ for all $i$. 
If $\sum_j\al_j=1$ then $k=1$ and
$e_{(\al_1,\dots,\al_k)}(a_1,\dots,a_k)=e_1(a_1)$. 
Suppose the claim true for all
$e_{(\bt_1,\dots,\bt_h)}(b_1,\dots,b_h)$ with 
$b_1,\dots,b_h\in A_R(m)$ and $\sum_i\bt_i < \sum_j\al_j$.
Let $k, a_1,\dots ,a_k, \al$ be as in the statement, then we have by
Theorem 2.1
$$e_{\al_1}(a_1)e_{(\al_2,\dots,\al_k)}(a_2,\dots,a_k)=$$
$$=e_{(\al_1,\dots,\al_k)}(a_1,\dots,a_k)+\sum
e_{\gamma}(a_1,\dots,a_k,a_1a_2,\dots,a_1a_k),$$ 
where
$$\gamma=(\gamma_{10},\gamma_{01},\dots,\gamma_{0h},\gamma_{11},\gamma_{12},\dots,\gamma_{1h})$$
with $h=k-1$,
$\sum_{j=0}^{h} \gamma_{1j}=\al_1$ with 
$\sum_{j=1}^h \gamma_{1j}>0$, and 
$\gamma_{0j}+\gamma_{1j}=\al_j$ for $j=1,\dots,h$. 
Thus 
$$\gamma_{10}+\gamma_{01}+\dots+\gamma_{0h}+\gamma_{11}+\dots+\gamma_{1h}
=\sum_j\al_j - \sum_{j=1}^h \gamma_{1j}<\sum_j\al_j.$$ 
Hence  
$$e_{(\al_1,\dots,\al_k)}(a_1,\dots,a_k)=$$ 
$$e_{\al_1}(a_1)e_{(\al_2,\dots,\al_k)}(a_2,\dots,a_k)
-\sum e_{\gamma}(a_1,\dots,a_k,a_1a_2,a_1a_3,\dots,a_1a_k),$$
where $\sum_{r,s}\gm_{rs}<\sum_j\al_j$. So the claim follows by
induction hypothesis.
\enddemo

\example{Example 2.4}
Consider $e_{(2,1)}(a,b)$ in $A_R(3,m)$ as in Example 1.2, then 
$$e_{(2,1)}(a,b)=e_2(a)e_1(b)- e_{(1,1)}(a,ab)
=e_2(a)e_1(b)-e_1(a)e_1(ab)+e_1(a^2b).$$
\endexample

\noindent
We now recall some basic facts about classical symmetric functions,
for further reading on this topic see {\cite{M}}.  

\noindent
We have another distinguished kind of functions in $\Lambda_R$ beside 
the elementary symmetric ones: the {\it{power sums}}. 

\noindent
For any
$r\in \NN$ the $r$-th power sum is 
$$p_r:=\sum_{i\geq 1}x_i^r.$$ 

\noindent
Let $g\in \Lambda_R$, set 
$g \cdot p_r=g(x_1^r,x_2^r,\dots,x_k^r,\dots)$, this is again a symmetric
function. 
Since the $e_i$ generate $\Lambda_R$ we have that $g\cdot p_r$ can be
expressed as a polynomial in the $e_i$. In particular,  
$$P_{h,k}:=e_h\cdot p_k$$ 
is a polynomial in the $e_i$.

\proclaim{Proposition 2.5}
For all $f\in A_R(m)$, and $k,h\in \NN$, $e_h(f^k)$ belongs
to the subring of $A_R(n,m)^{S_n}$ generated by the $e_j(f)$.
\endproclaim
\demo{Proof} 
Let $f\in A_R(m)$ and consider $e_h(f^k)\in A_R(n,m)^{S_n}$, we have
(see Introduction) 
$$e_h(f^k)=e_h(f(1)^k,\dots,f(n)^k)=
P_{h,k}(e_1(f(1),\dots,f(n)),\dots,e_n(f(1),\dots,f(n)))$$
and the result is proved.
\enddemo

\noindent
We are now ready to prove Theorem 1 stated in the introduction.

\demo{Proof of Theorem 1} Recall that a monomial $\mu\in \mn^+$ is
called {\it{primitive}} if it is not a power of another one and we
denote by $\frak M_m^+$ the set of primitive monomials. 
The elements $e_{\al}\in\Cal B_{n,m,R}$, that form a $R$-basis by
Prop.1.4, can be expressed as polynomials in $e_i(\mu)$ with
$i=1,\dots, n$ and $\mu \in \mm$, by Cor.2.3. If $\mu=\nu^k$ with
$\nu\in \frak M_m^+$, then $e_i(\mu)$ can be expressed as a polynomial
in the $e_j(\nu)$, by Prop.2.5. Since for all
$\mu\in\mm$ there exist $k\in\NN$ and $\nu\in\frak M_m^+$ such that
$\mu=\nu^k$, we have that $A(n,m)^{S_n}$ is generated as a commutative
ring by the $e_j(\nu)$, where $\nu\in \frak M_m^+$ and $j=1,\dots ,n$. 

\noindent
The theorem then follows by the following result due to Fleischmann
{\cite{F}}: the ring $A_R(n,m)^{S_n}$ is generated by elements of total
degree $\ell\leq n(m-1)$, for any commutative ring $R$, with
sharp bound if $n=p^s$ a power of a prime and $R=\ZZ$ or 
$p\cdot 1_R=0$. 
If $R\supset \QQ$ then the result follows from Newton's Formulas and
a well--known result of H.Weyl (see {\cite{G}},{\cite{W}}).\qed 
\enddemo

\head 3. Relations \endhead
\noindent
We write a generating series for the orbits of monomials
$$G(t):=\prod_{i=1}^n(1+\sum_{\mn^+}t_{\mu}\mu(i))
=\sum_{\al,\mid\al\mid\leq n}t^{\al}e_{\al},\eqno(3.1)$$
where $\al\in\sd {\NN} {\mn^+}$ and $t^{\al}e_{\al}(n)=0$ when $\al=0$.

\noindent
Recall the map $\pi_{n}:A_R(n,m)@>>> A_R(n-1,m)$ defined by
$$\pi_{n}(x_i(j))=\cases 0 \;&{\text{if}}\,\, j=n \\
x_i(j) \;&{\text{if}}\,\, j\leq n-1
\endcases \,\,\,\,\,\;\;\;\;\;\;\;{\text{ for all }}\, i.\eqno(3.2)
$$
Then we have of course that $\pi_{n}(G_n(t))=G_{n-1}(t)$, so that
$$\pi_{n}((e_{\al}))=\cases e_{\al} &\;{\text{if}}\; \mid \al \mid
																						\leq n\\            
																						0 &\;{\text{otherwise.}} 
\endcases
\eqno(3.3)$$ 
Thus, by Prop.1.4, for all $a\in\NN^m$ the restriction
$$\pi_{n,a} : A_R(n,m,a)\rightarrow A_R(n-1,m,a)\eqno(3.4)$$
is such that 
$$ \pi_{n,a}(A_R(n,m,a)^{S_n})=A_R(n-1,m,a)^{S_{n-1}}\eqno(3.5)$$
and then $(A_R(n,m,a)^{S_n},\pi_{n,a})$ is a projective sytem.

\noindent
For any $a \in \NN^m$ set 
$$A_R(\infty,m,a):=\lim_{\leftarrow } A_R(n,m,a)^{S_n}, \eqno(3.6)$$ 
where the projective limit is taken with respect to $n$ over the above
projective system and set 
$$\tilde{\pi}_{n,a}:A_R(\infty,m,a)@>>> A_R(n,m,a)^{S_n}\eqno(3.7)$$ 
for the natural projection.

\noindent
Set 
$$A_R(\infty,m):=\bigoplus_{a\in \NN^m} A_R(\infty,m,a)\eqno(3.8)$$
and
$$\tilde{\pi}_n:=\bigoplus_{a\in\NN^m}\tilde{\pi}_{n,a}.\eqno(3.9)$$

\noindent
Similarly to the classical case ($m=1$) and recalling (3.1), (3.3) we
make an abuse of notation and set 
$$e_{\al}:=\underset{\leftarrow }\to{\lim}\; e_{\al}(n),$$ 
for any $\al \in \sd {\NN}{\mm}$. In the same way we set 
$e_j(f):=\underset{\leftarrow }\to{\lim}\;e_j(f)$ with $j\in \NN$,
where $f\in A_R(m)^+$ is homogeneous of positive multidegree, so that
$j\;\partial(f)=a$.  

\proclaim{Proposition 3.1}Let $a\in \NN^m$. 
\roster
\item The $R$-module $\ker \tilde{\pi}_{n,a}$ is the linear span of 
$$\{e_{\al}\in A_R(\infty,m,a)\,:\,\mid \al \mid >n\}.$$
\item The $R$-module homomorphisms 
$\tilde{\pi}_{n,a}:A_R(\infty,m,a)@>>> A_R(n,m,a)^{S_n}$ are onto for
all $n\in\NN$ and $A_R(\infty,m,a) \cong A_R(n,m,a)^{S_n}$ for all 
$n \geq \mid a \mid.$ 
\item The $R$-module $A_R(\infty,m,a)$ is free with basis 
$$\{e_{\al} \, : \, \partial(e_{\al})=a \},$$
\item The $R$-module $A_R(\infty,m)$ is free with basis 
$$\{e_{\al} \, : \, \al\in\sd \NN {\mm} \}.$$
\endroster
\endproclaim

\demo{Proof}
\roster
\item By (3.3) and (3.5), for all $a\in\NN^m$, the following is a split
exact sequence of $R$-modules  
$$0 @>>> \ker {\pi}_{n,a} @>>> A(n,m,a)^{S_n} @>{\pi}_{n,a}>> 
A(n-1,m,a)^{S_{n-1}} @>>> 0,$$
and the claim follows.
\item If $\sum_{j=1}^m a_j < n $, then $\ker \tilde{\pi}_{n,a}= 0$,
indeed $$\partial(e_{\al})=\sum_{\mu \in \mm}
\al_{\mu}\;\partial(\mu)=a \Longrightarrow \mid \al \mid \leq
\sum_{j=1}^m a_j  < n.$$ 
Hence $A(h,m,a)^{S_h}\cong A(b,m,a)^{S_b}$ where $b:=\sum_{j=1}^m
a_j$, for all $h\geq \sum_{j=1}^m a_j$ and the claim follows by
(3.5). 
\item follows from (1) and (2).
\item follows from (3) and (3.8)\qed
\endroster

\enddemo

\noindent
\remark{Remark 3.2}
Notice that $A_R(m)^{\otimes n}\cong A_R(n,m)$ as multigraded
$S_n$-algebras by means of   
$$f_1\otimes\cdots \otimes f_n \leftrightarrow f_1(1)f_2(2)\cdots
f_n(n)\eqno(3.10)$$ 
for all $f_1,\dots ,f_n \in A_R(m)$.
Hence $A_R(n,m)^{S_n}\cong TS^n(A_R(m))$, where $TS^n(\;-\;)$ denotes
the symmetric tensors functor.  
Since $TS^n(A_R(m))\cong R\bigotimes TS^n(A_{\ZZ}(m))$ (see
{\cite{B}}),
we have 
$$A_R(n,m)^{S_n}\cong R\otimes A_{\ZZ}(n,m)^{S_n}\eqno(3.11)$$ 
for any commutative ring $R$.
\endremark

\noindent
We then work with $R=\ZZ$ and we suppress the ${\ZZ}$ subscript for
the sake of simplicity. 

\remark {Remark 3.3} The $\ZZ$-module $A(\infty,m)$ can be endowed
with a structure of $\NN^m$-graded ring such that the $\pi_n$ are
$\NN^m$-graded ring homomorphisms: the product $e_{\al}e_{\bt}$, where
$\al ,\bt \in \sd {\NN} {\mm}$, is defined by using the product formula
of Theorem 2.1 with no upper bound on $\mid\gamma\mid$, where $\gamma$
appears in the summation.  
\endremark

\proclaim{Proposition 3.4}
Consider the free polynomial ring 
$$C(m):=\bigoplus_{a\in \NN^m}C(m,a)
:=\ZZ[e_{i,\mu}]_{i\in \NN,\mu \in \frak M_m^+}$$ 
with multidegree given by $\partial(e_{i,\mu})=\partial(\mu)i$. 

\noindent
Then the multigraded ring homomorphism
$$\sigma_m:\ZZ[e_{i,\mu}]_{i\in \NN,\mu \in \frak M_m^+}@>>>A(\infty,m)$$
given by
$$\sigma_m: e_{i,\mu} \mapsto e_i(\mu), \,\,{\text{for all}}\,\,
i\in \NN,\mu \in \frak M_m^+$$
is an isomorphism, i.e. $A(\infty,m)$ is freely generated as a
commutative ring by the $e_i(\mu)$, where 
$i\in \NN$ and $\mu \in \frak M_m^+.$ 
\endproclaim
\demo{Proof}
Since we defined the product in $A(\infty,m)$ as in Theorem 2.1, it is
easy to verify, repeating the reasoning of the previous section,
that $A(\infty,m)$ is generated as a commutative ring by the $e_i(\mu)$,
where $i\in \NN,\mu \in \frak M_m^+.$ 
Hence $\sigma_m$ is onto for all $m\in\NN$.

\noindent
Let $a\in \NN^m$ and consider the restriction 
$\sigma_{m,a}:C(m,a)@>>>A(\infty,m,a)$.
It is onto as we have just seen. A $\ZZ$-basis of $C(m,a)$ is 
$$\{\prod_{i\in\NN,k\in\NN,\mu \in \frak M_m^+}e_{i,\mu}\,:\,
\sum_{i\in\NN,k\in\NN,\mu \in \frak M_m^+}i\;k\;\partial(\mu)\;=\;a\}.$$
On the other hand, a $\ZZ$-basis of $A(\infty,m,a)$ is 
$$\{e_{\al}\,:\,
\sum_{\al_{\mu}\in \NN,\mu \in \mm}\al_{\mu}\;\partial(\mu)\;=\;a\}.$$

\noindent
Let $\mu\in\mm$, then there are an unique $k\in\NN$ and an unique
$\nu\in\frak M_m^+$ such that $\mu=\nu^k$. Hence  
$$\sum_{\al_{\mu}\in \NN,\mu \in \mm}\al_{\mu}\;\partial(\mu) = 
\sum_{k\in\NN, \al_{\mu}\in \NN, \nu \in \frak M_m^+}\al_{\mu}\;k\;\partial(\nu),$$
so that $C(m,a)$ and $A(\infty,m,a)$ have the same (finite) $\ZZ$-rank
and thus are isomorphic via $\sigma_{m,a}$.\qed
\enddemo

\proclaim{Corollary 3.5}
Let $R\supset\QQ$ then $A_R(\infty,m)$ is a polynomial ring
freely generated by the $e_1(\mu)$, where $\mu\in\mm$. 
\endproclaim

\demo{Proof}
By Prop.3.4 and Theorem 1.\qed
\enddemo

\demo{Proof of Theorem 2} 
\roster
\item As before we set $R=\ZZ$ and the result follows by Remark
3.2, Prop.3.4. and Prop.3.1. 
\item
By Prop.3.1 the kernel of $$A(\infty,m) @>\tilde{\pi}_n>> A(n,m)^{S_n}$$ 
has basis $\{e_{\al}\:\mid \al \mid >n \}$. Let $V_k$ be the submodule
of $A(\infty,m)$ with basis $\{e_{\al}\:\mid \al \mid=k \}$. Let $A_k$
be the sub-$\ZZ$-module of $\QQ\otimes V_k$ generated by the $e_k(f)$
with $f\in A(m)^+$. Let $g:\QQ\otimes V_k @>>> \QQ$ be a linear form
identically zero on $A_k$. Then 
$$0=g(e_k(f))=g(e_k(\sum_{\mu\in \mm}\lambda_\mu \;\mu))
=(\sum_{\mid \al \mid=k}
(\prod_{\mu\in \mm}\lambda^{{\al}_{\mu}}_{\mu})\;g(e_{\al})),$$
for all $\sum_{\mu\in \mm}\lambda_{\mu}\;\mu \in A(m)^+$. Hence
$g(e_{\al})=0$ for all $e_{\al}$ with $\mid \al \mid=k$ ; thus $g=0$. 
If $R\supset\QQ$ the result then follows from Newton's formulas and
Cor.3.5.
\qed 
\endroster
\enddemo
\head Aknowledgement \endhead

\noindent
I would like to thank M.Brion, C.De Concini and C.Procesi, in
alphabetical order, for useful discussions.  
I would also like to thank the referee for its valuable suggestions.

\refstyle{C}

\enddocument